\documentclass{amsart}

\input{epsf.tex}
\usepackage{amsfonts,amssymb,verbatim,amsmath,amsthm,latexsym,textcomp,amscd}
\usepackage{latexsym,amsfonts,amssymb,epsfig,verbatim}
\usepackage{amsmath,amsthm,amssymb,latexsym,graphics,textcomp}
\usepackage{graphicx}
\usepackage{color}
\usepackage{url}

\input xy
\xyoption{all}

\begin{document}

\newtheorem{theorem}{Theorem}[section]
\newtheorem{prop}[theorem]{Proposition}
\newtheorem{lemma}[theorem]{Lemma}
\newtheorem{cor}[theorem]{Corollary}
\newtheorem{defn}[theorem]{Definition}
\newtheorem{conj}[theorem]{Conjecture}
\newtheorem{claim}[theorem]{Claim}
\newtheorem{defth}[theorem]{Definition-Theorem}

\newcommand{\boundary}{\partial}
\newcommand{\C}{{\mathbb C}}
\newcommand{\integers}{{\mathbb Z}}
\newcommand{\natls}{{\mathbb N}}
\newcommand{\ratls}{{\mathbb Q}}
\newcommand{\reals}{{\mathbb R}}
\newcommand{\proj}{{\mathbb P}}
\newcommand{\lhp}{{\mathbb L}}
\newcommand{\tube}{{\mathbb T}}
\newcommand{\cusp}{{\mathbb P}}
\newcommand\AAA{{\mathcal A}}
\newcommand\BB{{\mathcal B}}
\newcommand\CC{{\mathcal C}}
\newcommand\DD{{\mathcal D}}
\newcommand\EE{{\mathcal E}}
\newcommand\FF{{\mathcal F}}
\newcommand\GG{{\mathcal G}}
\newcommand\HH{{\mathcal H}}
\newcommand\II{{\mathcal I}}
\newcommand\JJ{{\mathcal J}}
\newcommand\KK{{\mathcal K}}
\newcommand\LL{{\mathcal L}}
\newcommand\MM{{\mathcal M}}
\newcommand\NN{{\mathcal N}}
\newcommand\OO{{\mathcal O}}
\newcommand\PP{{\mathcal P}}
\newcommand\QQ{{\mathcal Q}}
\newcommand\RR{{\mathcal R}}
\newcommand\SSS{{\mathcal S}}
\newcommand\TT{{\mathcal T}}
\newcommand\UU{{\mathcal U}}
\newcommand\VV{{\mathcal V}}
\newcommand\WW{{\mathcal W}}
\newcommand\XX{{\mathcal X}}
\newcommand\YY{{\mathcal Y}}
\newcommand\ZZ{{\mathcal Z}}
\newcommand\CH{{\CC\HH}}
\newcommand\PEY{{\PP\EE\YY}}
\newcommand\MF{{\MM\FF}}
\newcommand\RCT{{{\mathcal R}_{CT}}}
\newcommand\PMF{{\PP\kern-2pt\MM\FF}}
\newcommand\FL{{\FF\LL}}
\newcommand\PML{{\PP\kern-2pt\MM\LL}}
\newcommand\GL{{\GG\LL}}
\newcommand\Pol{{\mathcal P}}
\newcommand\half{{\textstyle{\frac12}}}
\newcommand\Half{{\frac12}}
\newcommand\Mod{\operatorname{Mod}}
\newcommand\Area{\operatorname{Area}}
\newcommand\ep{\epsilon}
\newcommand\hhat{\widehat}
\newcommand\Proj{{\mathbf P}}
\newcommand\U{{\mathbf U}}
 \newcommand\Hyp{{\mathbf H}}
\newcommand\D{{\mathbf D}}
\newcommand\Z{{\mathbb Z}}
\newcommand\R{{\mathbb R}}
\newcommand\Q{{\mathbb Q}}
\newcommand\E{{\mathbb E}}
\newcommand\til{\widetilde}
\newcommand\length{\operatorname{length}}
\newcommand\tr{\operatorname{tr}}
\newcommand\gesim{\succ}
\newcommand\lesim{\prec}
\newcommand\simle{\lesim}
\newcommand\simge{\gesim}
\newcommand{\simmult}{\asymp}
\newcommand{\simadd}{\mathrel{\overset{\text{\tiny $+$}}{\sim}}}
\newcommand{\ssm}{\setminus}
\newcommand{\diam}{\operatorname{diam}}
\newcommand{\pair}[1]{\langle #1\rangle}
\newcommand{\T}{{\mathbf T}}
\newcommand{\inj}{\operatorname{inj}}
\newcommand{\pleat}{\operatorname{\mathbf{pleat}}}
\newcommand{\short}{\operatorname{\mathbf{short}}}
\newcommand{\vertices}{\operatorname{vert}}
\newcommand{\collar}{\operatorname{\mathbf{collar}}}
\newcommand{\bcollar}{\operatorname{\overline{\mathbf{collar}}}}
\newcommand{\I}{{\mathbf I}}
\newcommand{\tprec}{\prec_t}
\newcommand{\fprec}{\prec_f}
\newcommand{\bprec}{\prec_b}
\newcommand{\pprec}{\prec_p}
\newcommand{\ppreceq}{\preceq_p}
\newcommand{\sprec}{\prec_s}
\newcommand{\cpreceq}{\preceq_c}
\newcommand{\cprec}{\prec_c}
\newcommand{\topprec}{\prec_{\rm top}}
\newcommand{\Topprec}{\prec_{\rm TOP}}
\newcommand{\fsub}{\mathrel{\scriptstyle\searrow}}
\newcommand{\bsub}{\mathrel{\scriptstyle\swarrow}}
\newcommand{\fsubd}{\mathrel{{\scriptstyle\searrow}\kern-1ex^d\kern0.5ex}}
\newcommand{\bsubd}{\mathrel{{\scriptstyle\swarrow}\kern-1.6ex^d\kern0.8ex}}
\newcommand{\fsubeq}{\mathrel{\raise-.7ex\hbox{$\overset{\searrow}{=}$}}}
\newcommand{\bsubeq}{\mathrel{\raise-.7ex\hbox{$\overset{\swarrow}{=}$}}}
\newcommand{\tw}{\operatorname{tw}}
\newcommand{\base}{\operatorname{base}}
\newcommand{\trans}{\operatorname{trans}}
\newcommand{\rest}{|_}
\newcommand{\bbar}{\overline}
\newcommand{\UML}{\operatorname{\UU\MM\LL}}
\newcommand{\EL}{\mathcal{EL}}
\newcommand{\tsum}{\sideset{}{'}\sum}
\newcommand{\tsh}[1]{\left\{\kern-.9ex\left\{#1\right\}\kern-.9ex\right\}}
\newcommand{\Tsh}[2]{\tsh{#2}_{#1}}
\newcommand{\qeq}{\mathrel{\approx}}
\newcommand{\Qeq}[1]{\mathrel{\approx_{#1}}}
\newcommand{\qle}{\lesssim}
\newcommand{\Qle}[1]{\mathrel{\lesssim_{#1}}}
\newcommand{\simp}{\operatorname{simp}}
\newcommand{\vsucc}{\operatorname{succ}}
\newcommand{\vpred}{\operatorname{pred}}
\newcommand\fhalf[1]{\overrightarrow {#1}}
\newcommand\bhalf[1]{\overleftarrow {#1}}
\newcommand\sleft{_{\text{left}}}
\newcommand\sright{_{\text{right}}}
\newcommand\sbtop{_{\text{top}}}
\newcommand\sbot{_{\text{bot}}}
\newcommand\sll{_{\mathbf l}}
\newcommand\srr{_{\mathbf r}}
\newcommand\geod{\operatorname{\mathbf g}}
\newcommand\mtorus[1]{\boundary U(#1)}
\newcommand\A{\mathbf A}
\newcommand\Aleft[1]{\A\sleft(#1)}
\newcommand\Aright[1]{\A\sright(#1)}
\newcommand\Atop[1]{\A\sbtop(#1)}
\newcommand\Abot[1]{\A\sbot(#1)}
\newcommand\boundvert{{\boundary_{||}}}
\newcommand\storus[1]{U(#1)}
\newcommand\Momega{\omega_M}
\newcommand\nomega{\omega_\nu}
\newcommand\twist{\operatorname{tw}}
\newcommand\modl{M_\nu}
\newcommand\MT{{\mathbb T}}
\newcommand\Teich{{\mathcal T}}
\renewcommand{\Re}{\operatorname{Re}}
\renewcommand{\Im}{\operatorname{Im}}

\title{Addendum to Ending Laminations and Cannon-Thurston Maps: Parabolics}

\author{Shubhabrata Das}

\author{Mahan Mj}
\address{RKM Vivekananda University, Belur Math, WB-711 202, India}

\thanks{This article is part of SD's PhD thesis  written under the supervision of MM. Research of the second author is supported in part by a CEFIPRA research grant.}

\date{}

\begin{abstract}
In earlier work, we had shown that Cannon-Thurston
 maps exist for Kleinian punctured surface groups without accidental parabolics. 
In this note we prove that pre-images of points are precisely end-points of leaves of the ending
lamination whenever the Cannon-Thurston map is not one-to-one.  This extends earlier work done for closed surface groups.  
\end{abstract}

\maketitle

\begin{center}
AMS subject classification =   57M50
\end{center}

\section{Ingredients} In \cite{mahan-split} we proved the existence of Cannon-Thurston maps for surface groups without accidental parabolics.
For closed surface groups, we described the structure of these maps in terms of ending lainations in \cite{mahan-elct}. In this note
we extend the results of  \cite{mahan-elct} to punctured surfaces.

\noindent {\bf 1) Equivalence Relations on ${\mathbb{S}}^1$:} \\ Suppose that a group $G$ acts on $S^1$ preserving a closed
 equivalenve relation $\mathcal R$. An example $\RR_\LL$ of such a relation comes from a lamination $\LL$, where two points on $S^1$
are equivalent if they are end-points of a leaf of $\LL$. The equivalence relation $\RR_\LL$ is obtained as the transitive closure of this relation.

\begin{defn} \cite{bowditch-ct} Two disjoint subsets, $P, Q \subset S^1$ are linked if there exist linked pairs, $\{x, y \} \subset P$ and $\{ z, w \} \subset Q$.
$\RR$ is unlinked if distinct equivalence classes are unlinked. \end{defn}

The  following Lemma due to Bowditch give us a way of recognising relation coming from laminations.

\begin{lemma} (Lemma 9.2 of \cite{bowditch-ct}) Let $\RR$ be a non-empty closed unlinked $G$-invariant equivalence relation on
$S^1$. Suppose that no pair of fixed points of any loxodromics are identified under $\RR$. Then
there is a unique complete perfect lamination, $\LL$, on $S$ such that $\RR = \RR_\LL$. \label{unlink} \end{lemma}

\noindent {\bf 2) Partial Electrocution}

Let $Y$ be the convex core of a simply (resp. doubly) degenerate 3-manifold of the form $S \times J$, where $J = [0, \infty )$ (resp.
$\mathbb R$).
Let $\mathcal{B}$ denote the equivariant collection of horoballs in $\til Y$ covering the cusps of $Y$. Let $X$ denote
$Y$ minus the interior of the horoballs in $\mathcal{B}$. Let 
$\mathcal{H}$ denote the collection of boundary horospheres.Then each
$H \in \mathcal{H}$ with the induced metric is isometric to a Euclidean
product $\mathbb{R} \times J$. 
Partially electrocute  each 
$H \in \HH$ by giving it the product of the zero metric (in the $\mathbb R$-direction) with the Euclidean metric
 (in the $J$-direction), The resulting space is quasi-isometric to what one would get by gluing
to each $H$ the mapping cylinder of the projection of $H$ onto the $J$-factor. Let $\JJ$ denote the collection of copies of $J$ obtained in this latter
construction and let $(\PEY , d_{pel})$ denote the resulting partially electrocuted space. (See \cite{mahan-pal} for a more general discussion.)
We have the following basic Lemma.

\begin{lemma} \cite{mahan-pal}
$(\PEY ,d_{pel})$ is a hyperbolic metric space and the sets $J_\alpha \in \JJ$
are uniformly quasiconvex.
\label{pel}
\end{lemma}

\noindent {\bf 3) Split Geometry and Ladders}

As pointed out in \cite{mahan-split} there exist a sequence of surfaces exiting the end(s) of $Y$, such that after removing the 
cusps and partially electrocuting them, and subsequently electrocuting split blocks,
we obtain a model of split geometry. This was utilized in \cite{mahan-split} to obtain the structure
of Cannon-Thurston maps for closed surface groups without accidental parabolics. The construction of split geometry recalled in the main body of \cite{mahan-split}
goes through for punctured surfaces with the pro viso that we first partially electrocute the space. We refer to \cite{mahan-split} for details.
Let $S_i$ be the sequence of truncated surfaces (i.e. surfaces minus cusps)
 exiting the end and forming boundaries of the blocks in the model of split geometry. Note that after partial
electrocution, the induced metric on the boundary horocycles of each $\til{S_i}$ is  the zero metric.
Equivalently each horocycle $HC$ is coned off to the corrsponding point of a $J_\alpha \in \JJ$. Let $(\til{S_i},d_{el})$ denote the resulting
electric metrics.

From a geodesic $\lambda = \lambda_0 \subset (\tilde{S_0},d_{el}) 
\subset \PEY$ we constructed in \cite{mahan-split} a `hyperbolic ladder'
${\LL}_\lambda$ such that $\lambda_i  =  {\LL}_\lambda \cap \tilde{S_i}$ is an
electro-ambient quasigeodesic
in the (path) electric metric on $ \tilde{S_i}$ induced by the graph metric $d_G$ on $\PEY$. 
 
We also constructed a large-scale retract $\Pi_\lambda : \PEY \rightarrow {\LL}_\lambda$ such that
the restriction $\pi_i$ of $\Pi_\lambda$ to $\tilde{S} \times \{ i \}$ is, roughly speaking,
a nearest-point retract of $\tilde{S} \times \{ i \}$ onto $\lambda_i$ in the (path) electric 
metric on $ \tilde{S_i}$.

We have the following basic theorem from \cite{mahan-split}.

\begin{theorem} \cite{mahan-split}
There exists $C > 0$ such that 
for any geodesic $\lambda = \lambda_0 \subset \widetilde{S} \times \{
0 \} \subset \widetilde{B_0}$, the retraction $\Pi_\lambda :
\PEY \rightarrow {\LL}_\lambda $ satisfies: \\

Then $d_{G}( \Pi_{\lambda } (x), \Pi_{\lambda } (y)) \leq C
d_{G}(x,y) + 
C$.
\label{retract}
\end{theorem}

\noindent{\bf  4) qi Rays}\\
We also have the following from \cite{mahan-split}.

\begin{lemma}
There exists $C \geq 0$ such that for $x_{i} \in \lambda_{i}$ there
exists
$x_{i-1} \in \lambda_{i-1}$ with $d_G (x_{i}, x_{i-1}) \leq C$. Similarly 
there
exists
$x_{i+1} \in \lambda_{i+1}$ with $d_G (x_{i}, x_{i+1}) \leq C$. Hence, for all $n$ and $x \in \lambda_n$,
 there exists a $C$-quasigeodesic ray $r$ such that $r(i) \in \lambda_i \subset {\LL}_\lambda$ for all $i$
and $r(n) = x$. 
\label{dGqgeod}
\end{lemma}

Hence, given $p \in \lambda_{i}$ the sequence of points $ x_{n}, n \in \mathbb{N} \cup \{ 0 \}$ (for simply degenerate groups)
or $n \in \integers$ (for totally degenerate groups) with $x_i = p$ gives by Lemma \ref{dGqgeod} above, a quasigeodesic in
the $d_G$-metric. Such quasigeodesics shall be referred to as {\em $d_G$-quasigeodesic rays}. Recall the following Proposition from \cite{mahan-elct}
in this context.

\begin{prop}
Let $\mu$, $\lambda$ be two bi-infinite geodesics on $\til{S}$
such that $\mu \cap \lambda \neq \emptyset$. Then ${\LL}_\lambda
\cap \LL_\mu = r$ is a quasigeodesic ray in $(\PEY, d_G)$..
\label{twin}
\end{prop}

\smallskip

\noindent {\bf 5) Easy Direction: Ideal points  are identified by  Cannon-Thurston Maps }
The easy direction of the main Theorem \ref{main} of this appendix is the same as that in \cite{mahan-elct}.

\begin{prop}
 Let $u, v$ be either ideal 
end-points of a leaf of a lamination, or ideal boundary points of a
complementary ideal polygon. Then
 $\partial i(u) = \partial i(v)$.
\label{ptpreimage}
\end{prop}

As in \cite{mahan-elct}
a {\bf CT leaf} $\lambda_{CT}$ will be a bi-infinite geodesic whose end-points are identified by $\partial i$. \\ An {\bf EL leaf} $\lambda_{EL}$ is a bi-infinite geodesic whose end-points
are  ideal boundary points of
either a leaf of the ending lamination, or a complementary ideal polygon. $\Lambda$ will denote the ending lamination for a simply degenerate group.

It remains to show that  \\
$\bullet$ {\bf A {\it CT leaf} is an {\it EL leaf}.}

\smallskip

Let $\RCT$ denote the equivalence relation on $S^1$ induced by the Cannon-Thurston map for the punctured surface (existence proven in \cite{mahan-split}.
By Proposition \ref{ptpreimage} pairs of end-points of leaves of $\Lambda$ are contained in $\RCT$. Hence, for simply degenerate groups,
 it suffices to show that 
$\RCT$ is induced by a lamination (since no other lamination can properly contain $\Lambda$). By Lemma \ref{unlink}
it suffices to show that $\RCT$ is unlinked. This is the content of the next section.

\section{$\RCT$ is unlinked}

We adapt Corollary 2.7 of \cite{mahan-elct} to the present context.

\begin{prop} \cite{mahan-elct} {\bf CT leaves have infinite diameter} \\
Let $\lambda_+ ( \subset \lambda \subset \tilde{S} \times \{ 0 \} = \tilde{S}) = {\mathbb H}^2$ be a semi-infinite geodesic (in the {\em hyperbolic metric} on $\tilde{S}$) contained in a {\it CT leaf} $\lambda$. Further suppose that $\lambda_+$ does not have a parabolic as its limit point in $\partial {\mathbb H}^2$.
Then $dia_G (\lambda_+)$ is infinite, where $dia_G$ denotes diameter in the graph metric restricted to 
$\til{S}$.
\label{ctleafinf}
\end{prop}

\begin{prop} Let $i: {\til S} \rightarrow \til{M}$ be an inclusion of the universal cover of a punctured surface
into the universal cover of the convex core $M$ of a simply degenerate 3-manifold. Let $\partial i$ be the associated Cannon-Thurston map.
If $\lambda$ is a CT-leaf in $\til S$, $\LL_\lambda$ the corrsponding ladder,  and $r=r(n) \subset \LL_\lambda$ a qi  ray, then there exists
$z\in\partial \til{M}$ such that $r(n) \rightarrow z = \partial i(\lambda_{-\infty }) =
\partial i(\lambda_{\infty })$ as $n \rightarrow \infty$. 
\label{qgeodasymp}
\end{prop}

\noindent {\bf Proof:} We first observe that both end-points $\lambda_{-\infty }, \lambda_{\infty }$ of the CT leaf $\lambda$ cannot be parabolics.
For then they would have to be base points of {\it different} horoballs in $\til M$ as they correspond to different lifts of the cusp(s) of $M$. 

\smallskip

\noindent {\bf Case 1:} Both $\lambda_{-\infty }, \lambda_{\infty }$ are non-parabolic. \\
The proof  of Proposition 2.13 of \cite{mahan-elct} goes through in this context mutatis mutandis. 

\smallskip

\noindent {\bf Case 2:} Exactly one of $\lambda_{-\infty }, \lambda_{\infty }$  is a parabolic. \\
Without loss of generality assume that  $\lambda_{-\infty }$ is a parabolic. 
Let $B$ be the horoball in $\til M$ based at $w = \partial i (\lambda_{-\infty} )$ and let $H$ be the horosphere boundary of $B$. Let $o$ be
the point of intersection of $\lambda$ with $H$. Choose a sequence of points $a_n, b_n \in \lambda$
such that $a_n \rightarrow \lambda_-\infty $
and $b_n \rightarrow \lambda_\infty $. Let $\overline{a_nb_n}$ be the geodesic in $\til M$ joining $a_n, b_n$. Then by the existence
of Cannon-Thurston maps for $i: \til{S} \rightarrow \til{M}$  it follows easily (see Lemma 2.1 of \cite{mitra-trees} for instance) that there exists
a function $M(n) \rightarrow \infty$ as $n \rightarrow \infty$ such that $\overline{a_nb_n}$ lies outside $B_{M(n)} (o) \subset \til{M}$.
Hence, if $q_n = \overline{a_nb_n} \cap H$ then $d(q_n, o) \geq M(n)$ and the geodesic subsegment $\overline{q_nb_n}$ lies outside $B_{M(n)} (o) \subset \til{M}$.

Let $N = \til{M} \setminus \bigcup_\alpha B_\alpha$ be the complement of open horoballs and $d_G$ be the graph metric on $N$ obtained after first
partially electrocuting horospheres. Let $(q_n,b_n)$ be the geodesic joining $q_n, b_n$ in $(N,d_G)$. Then by weak relative hyperbolicity of $N$,
(see Lemma 2.1 of \cite{mahan-split}) $(q_n,b_n)$ and $\overline{q_nb_n}$ lie in a bounded neighborhood of each other in $(N,d_G)$. 

By Proposition \ref{ctleafinf} $d_G(o,b_n)\rightarrow \infty$ as $n\rightarrow \infty$. Also, $d_G(o, q_n)$ is the number of vertical
blocks between $o$ and $q_n$ and hence $d_G(o, q_n) \rightarrow \infty$. But $a_n \rightarrow \lambda_-\infty$ implies
$q_n \rightarrow \infty$. Hence $d_G(o, q_n) \rightarrow \infty$ as $a_n \rightarrow \lambda_-\infty$. Finally (see for instance the argument in Sections
6.3, 6.4 of \cite{mahan-split}) there exists
a function $M_1(n) \rightarrow \infty$ as $n \rightarrow \infty$ such that $(q_n,b_n)$ lies outside $B_{M_1(n)} (o) \subset (N,d_G)$.

Now recall that $\Pi_\lambda : N \rightarrow \LL_\lambda$ is a coarse Lipschitz retract by Theorem \ref{retract}. Hence $\Pi_\lambda [(q_n,b_n)]
\subset \LL_\lambda$ uniform quasigeodesic in $(N,d_G)$.

Further, since $q_n \in H$ and since $\Pi_\lambda$ essentially fixes the horosphere $H$, it follows that $d_G( \Pi_\lambda (q_n), q_n) \leq 1$.
Also $\Pi_\lambda (b_n) = b_n$. Therefore there exists
a function $M_2(n) \rightarrow \infty$ as $n \rightarrow \infty$ such that $\Pi_\lambda [(q_n,b_n)]$ lies outside $B_{M_2(n)} (o) \subset (N,d_G)$.

Next, since $H \cap \LL_\lambda$ and $b_n$ lie on opposite sides of the  qi  ray $r=r(n) \subset \LL_\lambda$ it follows that there exists
$z_n \in (q_n,b_n)$ such that $d_G(z_n,r)$ is uniformly bounded.

Also there exists $t_n \in \overline{q_nb_n}$ such that  $d_G(z_n,t_n)$ and hence $d_G(t_n,r)$ is uniformly bounded.

Since $t_n \in \overline{q_nb_n}$ it follows that $t_n \rightarrow \partial i (\lambda_{-\infty} ) = \partial i (\lambda_{\infty} )$.
Since $d_G(t_n,r)$ is uniformly bounded, there exists $s_n \in r$ such that
$d_G(t_n,s_n)$ is uniformly bounded and therefore $t_n, s_n$ are separated by a uniformly bounded number of split components. By uniform graph
quasiconvexity of split components (Definition-Theorem 1.6 of \cite{mahan-elct}) it follows that $s_n \rightarrow \partial i (\lambda_{-\infty} ) = \partial i (\lambda_{\infty} )$. Finally if $r_{s_n}$ denotes the part of the ray $r$ `above' $s_n$, (i.e. $[s_n, \infty)$) then joining  points of $r_{s_n}$ in
successive blocks by hyperbolic geodescics we obtain a path $\sigma_n$ which contains a semi-infinite hyperbolic ray (the limit of hyperbolic geodesic
segments joining $s_n$ to points arbitrarily far along $r$) (by weak relative hyperbolicity of $N$,
-- Lemma 2.1 of \cite{mahan-split}). Hence  $r(n) \rightarrow  \partial i(\lambda_{-\infty }) =
\partial i(\lambda_{\infty })$ as $n \rightarrow \infty$.
$\Box$

\begin{theorem} \label{main} 
 Let $\partial i(a) = \partial i(b)$ for $a, b \in
S^1_\infty$ be two distinct points that are identified by the
Cannon-Thurston map corresponding to a simply degenerate  surface
group (without accidental parabolics). 
Then $a, b$ are either ideal 
end-points of a leaf of the ending lamination (in the sense of Thurston), 
or ideal boundary points of a
complementary ideal polygon. Conversely, if $a, b$ are either ideal 
end-points of a leaf of a lamination, or ideal boundary points of a
complementary ideal polygon, then $\partial i(a) = \partial i(b)$. \end{theorem}

Suppose $\lambda$ and $\mu$ are intersecting CT leaves, i.e.
$\partial i(\lambda_{-\infty }) =
\partial i(\lambda_{\infty })$  and $\partial i(\mu_{-\infty }) =
\partial i(\mu_{\infty })$. 

Consider the ladders ${\LL}_\lambda$ and ${\LL}_\mu$. Let $r(i) = \lambda_i \cap \mu_i$ be a 
quasigeodesic ray as per Proposition \ref{twin}. 
By Proposition \ref{qgeodasymp}, $r$ converges
to a point $z$ on $\partial \til{M}$ 
such that $z =  \partial i(\lambda_{-\infty }) =
\partial i(\lambda_{\infty })  = \partial i(\mu_{-\infty }) =      
\partial i(\mu_{\infty })$. Hence if $ \{ a, b \},  \{ c, d \} \in \RCT$, then either $ \{ a, b, c, d \} $ are all mutually related in $\RCT$, or 
$ \{ a, b \},  \{ c, d \}$ are unlinked. By Lemma \ref{unlink}, $\RCT$ is induced by a lamination $\Lambda_{CT}$. By Proposition \ref{ptpreimage}, 
the ending lamination $\Lambda_{EL}$ is contained in $\Lambda_{CT}$. Since $\Lambda_{EL}$ is filling and arational, it follows that
$\Lambda_{EL} = \Lambda_{CT}$.
$\Box$

\medskip

The 
modifications necessary to pass from the simply degenerate case to the totally degenerate case are exactly as in the last subsection of \cite{mahan-elct}.

\bibliography{append}

\begin{thebibliography}{Bow07}

\bibitem[Bow07]{bowditch-ct}
B.~H. Bowditch.
\newblock The {C}annon-{T}hurston map for punctured surface groups.
\newblock {\em Math. Z. 255}, pages 35--76, 2007.

\bibitem[Mit98]{mitra-trees}
Mahan Mitra.
\newblock Cannon-{T}hurston {M}aps for {T}rees of {H}yperbolic {M}etric
  {S}paces.
\newblock {\em Jour. Diff. Geom.48}, pages 135--164, 1998.

\bibitem[Mj06]{mahan-split}
Mahan Mj.
\newblock {Cannon-Thurston Maps for Surface Groups}.
\newblock {\em preprint, arXiv:math.GT/0607509}, 2006.

\bibitem[Mj07]{mahan-elct}
Mahan Mj.
\newblock {Ending Laminations and Cannon-Thurston Maps }.
\newblock {\em preprint, arXiv:math.GT/0702162}, 2007.

\bibitem[MP07]{mahan-pal}
Mahan Mj and Abhijit Pal.
\newblock {Relative Hyperbolicity, Trees of Spaces and Cannon-Thurston Maps}.
\newblock {\em arXiv:0708.3578, submitted to Geometriae Dedicata}, 2007.

\end{thebibliography}
\bibliographystyle{alpha}

\end{document}